\documentclass[11pt]{article}

\usepackage[margin=1in]{geometry}
\usepackage{amsmath,amssymb,amsthm,mathtools,mathrsfs}
\usepackage{enumitem}
\usepackage[colorlinks=true,linkcolor=blue,citecolor=blue,urlcolor=blue]{hyperref}
\mathtoolsset{showonlyrefs=true}

\newtheorem{theorem}{Theorem}
\newtheorem{lemma}{Lemma}
\newtheorem{proposition}{Proposition}

\DeclareMathOperator{\Vol}{Vol}
\DeclareMathOperator{\Tr}{tr}
\DeclareMathOperator{\Ric}{Ric}
\DeclareMathOperator{\Id}{Id}
\DeclareMathOperator{\SO}{SO}

\def\r{\mathbb{R}}
\def\rn{\mathbb{R}^n}
\def\s{\mathbb{S}^2}

\newcommand{\e}[1]{\begin{equation}#1\end{equation}}
\newcommand{\ald}[1]{\begin{aligned}#1\end{aligned}}

\newcommand{\ang}[1]{\left\langle#1\right\rangle}
\newcommand{\prn}[1]{\left(#1\right)}
\newcommand{\cur}[1]{\left\{#1\right\}}

\newcommand{\norm}[1]{\left\|#1\right\|}
\newcommand{\calH}{\mathcal{H}}
\newcommand{\calL}{\mathcal{L}}

\title{An Improved Lower Bound for the Three-Dimensional Blaschke--Lebesgue Problem from Spectral and Dual Perspectives}
\author{Akatsuki Nishioka\thanks{Department of Mathematical and Computing Science, School of Computing, Institute of Science Tokyo, Ookayama 2-12-1, Meguro-ku,
152-8550, Tokyo, Japan. Email: nishioka.a.2122@m.isct.ac.jp}}
\date{\today}

\begin{document}
\maketitle

\begin{abstract}
The Blaschke--Lebesgue problem asks for convex bodies of minimum volume among all convex bodies of prescribed constant width. In the plane, the minimizer is the Reuleaux triangle, whereas the corresponding three-dimensional problem remains open and is also known as Meissner's conjecture. In this paper, we establish the lower bound
$(4\pi/33)d^3 \simeq 0.380799\,d^3$ for the volume of any three-dimensional convex body of constant width d. This improves upon Chakerian's lower bound, approximately $0.364916\,d^3$, although it remains below the volume of the conjectured minimizers, Meissner's tetrahedra, whose volume is approximately $0.419860\,d^3$. The proof is based on a support-function formulation, spectral estimates via spherical harmonics, and Bochner's formula. We also show that the resulting lower bound can be interpreted as a Lagrange dual bound for the associated concave quadratic minimization problem. This dual viewpoint suggests possible routes toward sharper lower bounds.
\end{abstract}


\section{Introduction}

The classical Blaschke--Lebesgue theorem states that the Reuleaux triangle minimizes the area among all plane convex bodies of prescribed constant width. It was proved independently in 1914 by Lebesgue \cite{lebesgue1914} and in 1915 by Blaschke \cite{blaschke1915}. Since then, a number of alternative proofs have appeared \cite{bogosel25new,ghandehari96,harrell02,hynd24,malagoli09}.

The corresponding three-dimensional problem remains open. In 1934, Bonnesen and Fenchel \cite{bonnesen1934} conjectured that Meissner's tetrahedra\footnote{There are two types of Meissner's tetrahedra with the same volume. See monographs \cite{kawohl11,martini19} for details.} minimize the volume among three-dimensional convex bodies of constant width. This conjecture is also known as Meissner's conjecture.\footnote{This terminology is used, for example, in \cite{antunes22,bogosel24mixed,oudet13}. Kawohl and Weber \cite{kawohl11} note that they could not find a written record in which Meissner himself explicitly stated the conjecture, although he seems to have expected that Meissner's tetrahedra have minimal volume.} Since then, the problem has been studied from several analytical and numerical viewpoints \cite{anciaux11,antunes22,bayen12,bayen07,bogosel24mixed,bogosel26,chakerian66,oudet13}.

In this paper, we improve Chakerian's 1966 lower bound \cite{chakerian66} $(\pi/3)(3\sqrt{6}-7)d^3\simeq 0.364916\,d^3$ for the volume of three-dimensional convex bodies of constant width $d$. This bound is also mentioned in \cite{bogosel24mixed,kawohl11} and in \cite[Theorem 14.1.2]{martini19}. Our main result is the following.

\begin{theorem}\label{thm_main}
Let $K\subset\r^3$ be a convex body of constant width $d$. Then its volume satisfies
\e{
\Vol(K)\geq \frac{4\pi}{33}d^3\simeq 0.380799\,d^3.
}
\end{theorem}

This lower bound is still below the volume of Meissner's tetrahedra, which is approximately $0.419860\,d^3$. Note that the maximizer, the ball of diameter $d$, has volume $\pi d^3/6\simeq0.523599\,d^3$. Thus, there remains a large gap between the lower bound obtained here and the conjectured minimum.

The proof is based on a new support-function formulation, spectral estimates by spherical harmonics, and Bochner's formula in Riemannian geometry. It can also be interpreted as a Lagrange dual bound for the associated concave quadratic minimization problem. This viewpoint is useful because it connects the present argument with duality methods for nonconvex quadratic programming, which has been extensively studied in finite-dimensional optimization \cite{tuy16,xia11}, and suggests possible refinements of the lower bound.

The paper is organized as follows. Section \ref{sec_pre} recalls basic facts on support functions, bodies of constant width, spherical harmonics, and Bochner's formula. A new support-function formulation of the three-dimensional Blaschke--Lebesgue problem is also introduced. Section \ref{sec_spe} provides the spectral estimate that is used in the proof of the main theorem. Section \ref{sec_pro} proves Theorem \ref{thm_main}. Section \ref{sec_int} gives the Lagrange-dual interpretation of the lower bound. Finally, Section \ref{sec_dis} discusses possible directions for further improvement.

\subsection*{Notation}

$\rn$ denotes the $n$-dimensional Euclidean space with the standard inner product $\cdot$. $\mathbb{S}^{n-1}$ denotes the unit sphere of $\rn$. $g$ denotes the round metric of $\s$. $d\sigma$ denotes the surface measure of $\s$. $T_u \s$ denotes the tangent space of $\s$ at $u\in\s$. For a tangent vector $\xi_u\in T_u \s$ and a self-adjoint linear operator $A_u:T_u \s\to T_u \s$, we define $|\xi_u|=\sqrt{g(\xi_u,\xi_u)}$ and $|A_u|=\sqrt{\Tr(A_u A_u)}$. For a tangent vecter field $\xi$ and a tensor field $A$, we define $\|\xi\|^2_{L^2}=\int_{\s}|\xi|^2\,d\sigma$ and $\|A\|^2_{L^2}=\int_{\s}|A|^2\,d\sigma$.

\section{Preliminaries}\label{sec_pre}

\subsection{Support function}

Let $K\subset\r^3$ be a convex body (i.e., nonempty compact convex set). Its support function $p:\s\to\mathbb{R}$ is defined by
\e{
p(u)\coloneqq\sup_{x\in K}\,x\cdot u.
}
The width of $K$ in the direction $u$ is defined by $p(u)+p(-u)$. Thus, $K$ is said to have constant width $d>0$ if
\e{\label{eq:const}
p(u)+p(-u)=d,\quad\forall u\in\s.
}
In the following, we use the zero-mean version of the support function
\e{
h\coloneqq p-\frac{1}{4\pi}\int_{\s} p(u) \,d\sigma.
}
Integrating both sides of \eqref{eq:const} gives $\int_{\s} p(u) \,d\sigma=2\pi d$, and thus the constant-width condition \eqref{eq:const} becomes the following:
\e{\label{eq:odd}
h(-u)=-h(u),\quad\forall u\in\s.
}

A $C^2$ function $p:\s\to\r$ is the support function of a convex body if and only if the following tensor, whose eigenvalues are the principal radii of curvature, is nonnegative \cite{antunes22,howard06,schneider13}
\e{\label{eq:conv-supp}
\nabla^2 p(u)+p(u) \Id \succeq 0,\quad\forall u\in\s.
}
Here, $\nabla^2 p(u):T_u\s\to T_u\s$ is the self-adjoint linear map defined by $\nabla^2 p(u)X\coloneqq(\nabla_X\nabla p)(u)$ (the covariant derivative of the spherical gradient $\nabla p$ with respect to a vector field $X$ at $u$) and $\Id$ is the identity map on $T_u\s$. The above inequality is understood in the sense of quadratic forms on the tangent space $T_u\s$: for any $\xi\in T_u\s$, $g(\nabla^2_{\s}p(u)\xi,\xi)+p(u)g(\xi,\xi)\ge 0$. For general bodies of constant width, the support function is $C^{1,1}$ \cite{martini19}, and the Hessian is defined almost everywhere.

We define
\e{
A_h(u)\coloneqq\nabla^2 h(u)+h(u)\Id.
}
Then the convexity constraint \eqref{eq:conv-supp} is 
\e{
A_h(u)+\frac{d}{2} \Id\succeq0,\quad\forall u\in\s
}
Since $h$ is odd (so is $\nabla^2 h$ \cite[Lemma 3.5]{howard06}), we have\footnote{Although $A_h(u)$ and $A_h(-u)$ act on different tangent spaces $T_u\s$ and $T_{-u}\s$, these spaces are identified with the linear space $u^\perp\subset\r^3$, and these equalities and inequalities can be understood as equalities and inequalities of linear operators on $u^\perp$. See also \cite{howard06}.}
\e{
-A_h(u)+\frac{d}{2} \Id=A_h(-u)+\frac{d}{2} \Id\succeq0,\quad\forall u\in\s.
} 
Thus, $A_h$ satisfies
\e{\label{eq:matrix-box}
-\frac{d}{2}\Id\preceq A_h(u)\preceq \frac{d}{2}\Id,\quad\forall u\in\s.
}

\subsection{Spherical harmonics}\label{sec:sphe}

In the two-dimensional Blaschke--Lebesgue problem, the Fourier expansion of the support function plays an important role. In the three-dimensional case, the corresponding role is played by spherical harmonics. See \cite{dai13,martini19,schneider13} for details of spherical harmonics.

A spherical harmonic $Y_\ell$ of degree $\ell=0,1,2,\ldots$ on the sphere $\mathbb{S}^{n-1}$ is, by definition, the restriction to $\mathbb{S}^{n-1}$ of a homogeneous polynomial $f$ of degree $\ell$ on $\rn$ that satisfies $\Delta f=0$. When $n=2$, spherical harmonics reduce to the usual Fourier modes $1,\sin\theta,\cos\theta,\sin2\theta,\cos2\theta,\ldots$. From now on, we only use the case $n=3$. The spherical harmonics are eigenfunctions of the spherical Laplacian $\Delta$:
\e{
-\Delta Y_{\ell}=\lambda_\ell Y_{\ell},\qquad\lambda_\ell=\ell(\ell+1).
}
The space of square-integrable functions on $\s$ decomposes orthogonally into the linear subspaces $\calH_\ell\ (\ell=0,1,2,\ldots)$ of degree-$\ell$ spherical harmonics (including the zero function). In each space $\calH_\ell$, we choose an orthonormal basis $Y_{\ell ,1},\ldots,Y_{\ell,2\ell+1}$. Then every $f\in L^2(\s)$ admits the expansion
\e{
f=\sum_{\ell=0}^{\infty} \sum_{j=1}^{2\ell+1} \ang{f,Y_{\ell,j}}_{L^2(\s)} Y_{\ell,j},
}
where the convergence is in $L^2(\s)$.

Spherical harmonics have a definite parity:
\e{
Y_{\ell}(-u)=(-1)^\ell Y_{\ell}(u)\quad \forall u\in\s.
}
Consequently, the spherical harmonic expansion of an odd function contains only odd-degree terms.

\subsection{Volume formula}

We use the following formula for the volume of a convex body of constant width, proved in \cite{anciaux11}.

\begin{lemma}[{\cite[Theorem 2]{anciaux11}}]
Let $K$ be a convex body of constant width $d$ with zero-mean support function $h=p-(d/2)$. Then
\begin{equation}\label{eq:volume}
\Vol(K)=\frac{\pi d^3}{6}-\frac{d}{2}E(h),
\end{equation}
where
\begin{equation}\label{eq:E-def}
E(h)\coloneqq\int_{\s}\left(\frac12|\nabla h|^2-h^2\right)\,d\sigma.
\end{equation}
\end{lemma}

The Wirtinger inequality \cite[Proposition 1]{anciaux11} states that 
\e{
E(h)\ge0,\quad\forall h\in C^{1,1}(\s)\ \text{with}\ \int_{\s}h\,d\sigma=0
}
with equality if $h\in\calH_1$ (including $h\equiv 0$). Thus, on $\mathcal{X}$, $E$ is a nonnegative quadratic functional and hence convex.

For a support function, the degree-one spherical harmonics correspond to translations of the body. Since translations do not change either the width or the volume, we may impose $h\perp \calH_1$ in $L^2(\s)$.

We now rewrite the three-dimensional Blaschke--Lebesgue problem as an optimization problem whose variable is the zero-mean support function $h$. In the sequel, we work with smooth zero-mean support functions $h\in C^\infty$ in order to use the Hessian $\nabla^2h$ and Bochner's formula. The following approximation result justifies this restriction at the level of the infimum.

\begin{lemma}[{\cite[Theorem 8.5.1]{martini19} and \cite[Theorem 1.8.20]{schneider13}}]\label{lem:smooth}
Every convex body of constant width can be arbitrarily closely approximated in the Hausdorff metric by bodies of constant width with smooth support functions. Moreover, the volume functional is continuous with respect to the Hausdorff metric on the space of convex bodies.
\end{lemma}

Define the linear subspace 
\e{
\mathcal{X}\coloneqq\{h\in C^\infty(\s)\mid h\perp\calH_1,\ h(-u)=-h(u)\}
}
and the constrained set
\e{
\mathcal{C}\coloneqq\cur{h\in \mathcal{X}\,\middle|\,-\frac{d}{2} \Id\preceq A_h(u)\preceq \frac{d}{2} \Id,\ \forall u\in\s}.
}
Then the Blaschke--Lebesgue problem (volume minimization problem for bodies of constant width $d$) can be written, without changing its infimum, as the following concave quadratic minimization problem under convex constraints:
\e{\label{eq:p}
\inf_{h\in\mathcal{C}}\ \mathcal{V}(h)\coloneqq\frac{\pi d^3}{6}-\frac{d}{2}E(h).
}

The problem \eqref{eq:p} need not admit a minimizer in the smooth class. Nevertheless, by Lemma \ref{lem:smooth}, its infimum agrees with the infimum of the original three-dimensional Blaschke--Lebesgue problem.

To the best of the author's knowledge, the formulation \eqref{eq:p} is new in the literature. Harrell \cite{harrell02} introduced a closely related optimization formulation whose variable is the sum of the principal radii of curvature. In the two-dimensional case, this formulation is equivalent to the Blaschke--Lebesgue problem. In dimension three, however, it is a relaxation of the original problem, since some admissible functions need not correspond to convex bodies of constant width. In contrast, \eqref{eq:p} gives a direct support-function formulation of the three-dimensional Blaschke--Lebesgue problem at the level of the infimum.

\subsection{Bochner's formula}

We introduce Bochner's formula on $\s$ to evaluate the Hessian in $A_h$. See \cite{lee18,chow06} for details. 

\begin{lemma}[Bochner's formula on $\s$]\label{lem:bochner}
For every smooth function $f$ on $\s$,
\begin{equation}\label{eq:bochner-integrated}
        \int_{\s}|\nabla^2 f|^2\,d\sigma
        =
        \int_{\s}(\Delta f)^2\,d\sigma
        -
        \int_{\s}|\nabla f|^2\,d\sigma.
\end{equation}
\end{lemma}
\begin{proof}
Bochner's formula on a Riemannian manifold is
\e{
\frac12\Delta |\nabla f|^2
=
|\nabla^2 f|^2+g(\nabla f,\nabla\Delta f)+\Ric(\nabla f,\nabla f).
}
On the unit sphere $\s$, the Ricci curvature tensor is the standard metric: $\Ric=g$. Integrating over the closed manifold $S^2$ (which means $\int_{\s} \frac12\Delta |\nabla f|^2\,d\sigma=0$) and using integration by parts, we obtain (6).
\end{proof}

\section{A spectral bound by spherical harmonics and Bochner's formula}\label{sec_spe}

This section proves the inequality
\begin{equation}\label{eq:spectral-ineq}
        E(h)\leq \frac1{22}\norm{A_h}_{L^2}^2
\end{equation}
for $h\in\mathcal{X}$. We use spherical harmonics and Bochner's formula.

We first evaluate $E(h)$ and $A_h$ for spherical harmonics and then prove the desired inequality.

\begin{lemma}\label{lem:A-eigen}
Let $Y_\ell\in\calH_\ell$, so that $-\Delta Y_\ell=\lambda_\ell Y_\ell$.  Then
\begin{equation}\label{eq:E-eigen}
        E(Y_\ell)=\frac{\lambda_\ell-2}{2}\norm{Y_\ell}_{L^2}^2,
\end{equation}
and
\begin{equation}\label{eq:A-eigen}
        \norm{A_{Y_\ell}}_{L^2}^2=(\lambda_\ell-1)(\lambda_\ell-2)\norm{Y_\ell}_{L^2}^2.
\end{equation}
\end{lemma}

\begin{proof}
Since $-\Delta Y_\ell=\lambda_\ell Y_\ell$,
\e{
        \int_{\s}|\nabla Y_\ell|^2\,d\sigma
        =
        -\int_{\s}Y_\ell\Delta Y_\ell\,d\sigma
        =
        \lambda_\ell\norm{Y_\ell}_{L^2}^2,
}
which yields 
\e{
E(Y_\ell)=\int_{\s}\left(\frac12|\nabla Y_\ell|^2-Y_\ell^2\right)\,d\sigma=\left(\frac{\lambda_\ell}{2}-1\right)\norm{Y_\ell}_{L^2}^2.
}

Next,
\e{
        |A_{Y_\ell}|^2=|\nabla^2 {Y_\ell}+{Y_\ell}\Id|^2
        =|\nabla^2 {Y_\ell}|^2+2{Y_\ell}\Tr(\nabla^2 {Y_\ell})+{Y_\ell}^2|\Id|^2.
}
On $\s$, $\Tr(\nabla^2 {Y_\ell})=\Delta {Y_\ell}$ and $|\Id|^2=2$.  Hence
\e{\label{eq:ah}
\norm{A_{Y_\ell}}_{L^2}^2
=
\int_{\s}|\nabla^2 Y_\ell|^2\,d\sigma
+2\int_{\s}Y_\ell\Delta Y_\ell\,d\sigma
+2\int_{\s}Y_\ell^2\,d\sigma.   
}
By Lemma~\ref{lem:bochner},
\e{
        \int_{\s}|\nabla^2 Y_\ell|^2\,d\sigma
        =\int_{\s}(\Delta Y_\ell)^2\,d\sigma
        -
        \int_{\s}|\nabla Y_\ell|^2\,d\sigma   =\lambda_\ell(\lambda_\ell-1)\norm{Y_\ell}_{L^2}^2.
}
Therefore
\e{
        \norm{A_{Y_\ell}}_{L^2}^2
        =
        (\lambda_\ell(\lambda_\ell-1)-2\lambda_\ell+2)\norm{Y_\ell}_{L^2}^2
        =
        (\lambda_\ell-1)(\lambda_\ell-2)\norm{Y_\ell}_{L^2}^2.
}
\end{proof}

\begin{proposition}\label{prop:spectral}
If $h\in\mathcal{X}$, then
\e{
E(h)\leq \frac1{22}\norm{A_h}_{L^2}^2.
}
The equality holds if and only if $h\in\calH_3$.
\end{proposition}

\begin{proof}
We use the spherical harmonic expansion of $h\in\mathcal{X}$:
\e{
        h=\sum_{\ell\ \mathrm{odd},\ \ell\geq3} h_\ell,
        \qquad h_\ell\in\calH_\ell.
}
Since $h_{\ell_1}$ and $h_{\ell_2}$ are orthogonal and
\e{
\int_{\s}\nabla h_{\ell_1}\cdot \nabla h_{\ell_2}\,d\sigma = \lambda_{\ell_1}\int_{\s}h_{\ell_1}h_{\ell_2}\,d\sigma=0
}
if $\ell_1\neq\ell_2$, the only diagonal terms of $E(h)$ and $\norm{A_h}_{L^2}^2$ remain. Indeed, \eqref{eq:ah} with $h$ instead of $Y_\ell$ and Bochner's formula (Lemma~\ref{lem:bochner}) yield
\e{
\|A_{h}\|_{L^2}^2
=\int_{\s}(\Delta h)^2\,d\sigma-\int_{\s}|\nabla h|^2\,d\sigma+2\int_{\s} h\Delta h\,d\sigma+2\int_{\s}h^2\,d\sigma. 
}
All terms involving $h_{\ell_1}$ and $h_{\ell_2}$ with $\ell_1\neq\ell_2$ vanish due to $-\Delta h_\ell=\lambda_\ell h_\ell$ and the orthogonality of spherical harmonics.

Thus, by Lemma~\ref{lem:A-eigen}, we have
\e{
E(h)
=
\sum_{\ell\ \mathrm{odd},\ \ell\ge3}\frac{\lambda_\ell-2}{2}\norm{h_\ell}_{L^2}^2,
}
and
\e{
\norm{A_h}_{L^2}^2
=
\sum_{\ell\ \mathrm{odd},\ \ell\ge3}(\lambda_\ell-1)(\lambda_\ell-2)\norm{h_\ell}_{L^2}^2.
}
For each odd $\ell\ge3$,
\e{
\frac{(\lambda_\ell-2)/2}{(\lambda_\ell-1)(\lambda_\ell-2)}
=\frac1{2(\lambda_\ell-1)}\le \frac{1}{22}
}
since the smallest allowed $\lambda_\ell$ is $\lambda_3=3(3+1)=12$. Therefore, we obtain the inequality.

The equality holds if and only if $h_\ell\equiv0$ for all odd $\ell\ge5$,
which is equivalent to $h\in\mathcal H_3$.
\end{proof}

\section{Proof of Theorem \ref{thm_main}}\label{sec_pro}

\begin{proof}[Proof of Theorem \ref{thm_main}]

Let $h\in\mathcal{X}$ satisfy the constraint $-(d/2) \Id\preceq A_h(u)\preceq (d/2)\Id,\quad\forall u\in\s$. We have
\e{\label{eq:matrix-est}
|A_h(u)|^2\le 2\left(\frac{d}{2}\right)^2,\quad\forall u\in\s,
}
and
\e{
\|A_h\|_{L^2}^2=\int_{\s}|A_h|^2\,d\sigma\le 2\pi d^2.
}
Therefore, by Proposition \ref{prop:spectral} and the volume formula \eqref{eq:volume},
\e{
\mathcal{V}(h)=\frac{\pi d^3}{6}-\frac{d}{2}E(h)\ge \frac{\pi d^3}{6}-\frac{d}{2}\frac{1}{22}\|A_h\|_{L^2}^2 \ge \prn{\frac{1}{6}-\frac{1}{22}}\pi d^3 = \frac{4\pi }{33}d^3.
} Lemma \ref{lem:smooth} gives the same lower bound for all bodies of constant width.
\end{proof}

\section{A Lagrange dual interpretation}\label{sec_int}

\subsection{Derivation of the Lagrange dual problem}

We now explain in what sense the preceding estimate is the Lagrange dual bound.

We consider the Lagrange dual problem of \eqref{eq:p}. We first replace the convexity constraint \eqref{eq:matrix-box} by an equivalent quadratic constraint
\e{\label{eq:quad-const}
A_h^2-\prn{\frac{d}{2}}^2\Id\preceq0,\quad \forall u\in\s,
}
so that the Lagrange dual bound takes finite values. This technique is commonly used in finite-dimensional quadratic programming \cite{xia11}. 

The Lagrange multiplier for this constraint is a tensor field on $\s$ where $Z(u):T_u\s\to T_u\s$ is a self-adjoint linear operator satisfying
\e{
Z(u)\succeq0
}
for any $u\in\s$. We also assume the integrals below are finite for $Z$. In the sequel, we just write this constraint as $Z\succeq0$. We define the Lagrangian for Problem \eqref{eq:p} as follows:
\begin{equation}\label{eq:Lagrangian-Gamma}
\calL(h,Z)\coloneqq
\frac{\pi d^3}{6}-\frac{d}{2}E(h)+
\int_{\s}\Tr\prn{Z\prn{A_h^2-
\frac{d^2}{4}\Id}}\,d\sigma.
\end{equation}
For any $h\in\mathcal{X}$ satisfying \eqref{eq:quad-const} and any $Z\succeq0$,
\e{
\mathcal{V}(h)\ge\calL(h,Z),
}
and thus
\e{\label{eq:weak-duality}
\inf_{h\in\mathcal{C}} \mathcal{V}(h)\ge\sup_{Z\succeq0} \inf_{h\in\mathcal{X}}\calL(h,Z).
}
The right-hand side is called the Lagrange dual problem and gives a lower bound.

The functional $\inf_{h\in\mathcal{X}}\calL(h,Z)$ is finite if and only if the quadratic form
\e{\label{eq:dual-feasibility-general}
        h\longmapsto
        \int_{\s}\Tr(Z A_h^2)\,d\sigma-\frac{d}{2}E(h)
}
is nonnegative for all $h\in\mathcal{X}$. If this condition holds, then the infimum in the right-hand side of \eqref{eq:weak-duality} is attained at $h=0$, and thus the Lagrange dual problem can be written as
\e{\label{eq:dual}
\sup_{Z\succeq0}\ \frac{\pi d^3}{6}-\frac{d^2}{4}\int_{\s}\Tr Z\,d\sigma\quad\text{s.t.}\ 
\int_{\s}\Tr(ZA_h^2)\,d\sigma\geq \frac{d}{2}E(h),\quad\forall h\in\mathcal{X}.
}

\subsection{Reduction via rotational invariance}

We now show that, by the rotational invariance of the dual problem \eqref{eq:dual}, it is enough to consider multipliers of the form $Z=z\Id$ with $z\ge0$. This reduction allows us to solve the dual problem analytically.

For $\rho\in\SO(3)$, define the rotated function
\e{
h^\rho(u)\coloneqq h(\rho^{-1}u)
}
and the rotated multiplier
\e{
Z^\rho(u)\coloneqq D\rho_{\rho^{-1}u} Z(\rho^{-1}u) (D\rho_{\rho^{-1}u})^{-1},
}
where $D\rho_{\rho^{-1}u}:T_{\rho^{-1}u}\s\to T_u\s$ is the derivative of $\rho:\s\to\s$ at $\rho^{-1}u$. It is used so that $Z^\rho$ is a linear operator on $T_u\s$ instead of $T_{\rho^{-1}u}\s$. If $Z$ satisfies the constraints in \eqref{eq:dual}, then $Z^\rho$ also satisfies the constraints:
\e{\ald{
\int_{\s}\Tr(Z^\rho(u) A_h^2(u))\,d\sigma &= \int_{\s}\Tr(Z^\rho(\rho v) A_h^2(\rho v))\,d\sigma\\
&= \int_{\s}\Tr(D\rho_{v} Z(v) (D\rho_{v})^{-1} A_h^2(\rho v))\,d\sigma\\
&= \int_{\s}\Tr(Z(v) (D\rho_{v})^{-1} A_h^2(\rho v) D\rho_{v})\,d\sigma\\
&= \int_{\s}\Tr(Z(v)A_{h^{\rho^{-1}}}^2(v))\,d\sigma\\
&\geq \frac{d}{2}E(h^{\rho^{-1}})=\frac{d}{2}E(h),\quad\forall h\in\mathcal{X}.
}}
Here we used that $h\in\mathcal X$ implies $h^{\rho^{-1}}\in\mathcal X$ and that
$E(h^{\rho^{-1}})=E(h)$.
Moreover, $Z^\rho\succeq0$ and the objective value is unchanged:
\e{
\int_{\s}\Tr Z^\rho(u)\,d\sigma=\int_{\s}
\Tr(
D\rho_v Z(v)(D\rho_v)^{-1})\,d\sigma(v)=\int_{S^2}\Tr Z(v)\,d\sigma.
}

Let $\mu$ denote the Haar measure on $\mathrm{SO(3)}$ normalized so that $\mu(\mathrm{SO(3)})=1$; see, for
example, \cite{meckes19}.  We define the rotational average of $Z$ by
\e{
\overline{Z}(u)\coloneqq\int_{\SO(3)}Z^\rho(u)\,d\mu(\rho).
}
Since $Z^\rho(u)$ all belong to the fixed tangent space $T_u \s$, the integral is obtained by integrating each entry of the corresponding $2\times2$ matrix in any orthonormal basis of $T_u\s$. The constraints are convex and the objective is linear, so $\overline{Z}$ still satisfies the constraints and has the same objective value.

For $\tau\in\mathrm{SO(3)}$ such that $\tau u=u$, $D\tau_u:T_u\s\to T_u\s$ defines a rotation on $T_u\s$, i.e., $D\tau_u$ can be seen as an element of $\mathrm{SO(2)}$. $\overline{Z}(u)$ is invariant under any such $\tau$
\e{
D\tau_u \overline{Z}(u) (D\tau_u)^{-1}=\int_{\SO(3)} D\tau_u Z^\rho(u) (D\tau_u)^{-1}\,d\mu(\rho) = \int_{\SO(3)} Z^{\tau\rho}(u)\,d\mu(\rho) =\overline{Z}(u),
}
where we used the invariance of the Haar measure under the rotation $\tau$.
Thus, $\overline{Z}(u)$ must be a scalar multiple of the identity\footnote{Indeed, if a $2\times2$ symmetric matrix $\begin{pmatrix} a & b \\ b & c \end{pmatrix}$ is invariant under a $(\pi/2)$-rotation $\begin{pmatrix} 0 & 1 \\ -1 & 0 \end{pmatrix}$, then they commute $\begin{pmatrix} -b & a \\ -c & b \end{pmatrix}=\begin{pmatrix} a & b \\ b & c \end{pmatrix} \begin{pmatrix} 0 & 1 \\ -1 & 0 \end{pmatrix}=\begin{pmatrix} 0 & 1 \\ -1 & 0 \end{pmatrix}\begin{pmatrix} a & b \\ b & c \end{pmatrix}=\begin{pmatrix} b & c \\ -a & -b \end{pmatrix}$, and we obtain $b=0$ and $a=c$. We can apply this argument to $\overline{Z}(u):T_u\s\to T_u\s$ by fixing the orthonormal basis of $T_u\s$.}
\e{
\overline Z(u)=z(u)\Id
}
for some $z:\s\to\r$. Moreover, taking the trace, we have
\e{
2z(u) = \Tr \overline Z(u) = \int_{\SO(3)}\Tr Z^\rho(u)\,d\mu(\rho) = \int_{\SO(3)}\Tr Z(\rho^{-1}u)\,d\mu(\rho)=\frac{1}{4\pi}\int_{\s}\Tr Z(v)\,d\sigma,
}
since $v=\rho^{-1}u$ is distributed uniformly on $\s$.


Therefore, we can restrict the Lagrange multiplier $Z$ to be $z\Id$ for some scalar $z\ge0$, and the dual problem reduces to
\begin{equation}\label{eq:dual-scalar-z}
\sup_{z\ge0}\ \frac{\pi d^3}{6}-2\pi d^2 z\quad\text{s.t.}\ z\norm{A_h}_{L^2}^2\geq \frac{d}{2} E(h)\quad\forall h\in\mathcal{X}.
\end{equation}
By Proposition~\ref{prop:spectral}, the smallest admissible $z$ is
\e{
z_*=\frac{d}{44}.
}
Hence, the dual bound gives
\e{\label{eq:Gamma-bound}
\inf_{h\in\mathcal{C}}\mathcal{V}(h)\ge \frac{\pi d^3}{6}-
        2\pi d^2 z^* = \frac{4\pi}{33}d^3.
}

\section{Discussion}\label{sec_dis}

The present bound should not be expected to be sharp. Equality in the spectral estimate (Proposition \ref{prop:spectral}) would require $h\in\mathcal{H}_3$, whereas equality in the pointwise matrix estimate \eqref{eq:matrix-est} would require $|A_h(u)|^2=d^2/2$ almost everywhere. Equivalently, the principal radii of curvature would have to be either $0$ or $d$ almost everywhere. These two equalities are not expected to be compatible for the support function of a convex body of constant width.

For the present dual problem, the averaging argument shows that the dual multiplier reduces to a scalar multiplier $z\Id$. Thus, improving the bound requires going beyond this rotation-averaged relaxation. Possible directions include symmetry-breaking constraints and additional variables adapted to the structure of the Meissner tetrahedra. Such modified dual problems, however, need not admit explicit analytic solutions. Therefore, their values may have to be estimated by certified numerical methods. This suggests a possible direction for combining the present analytic approach with computer-assisted verification.

Problem \eqref{eq:p} has a structure similar to finite-dimensional concave quadratic minimization problems, such as the maximum cut problem of a graph \cite{xia11}. Duality and relaxation theories for finite-dimensional quadratic programming, such as the moment-SOS hierarchy \cite{henrion20}, may also be useful for further improving the lower bound.

\section*{Acknowledgements}

I am grateful to Professor Dorin Bucur for introducing me to \cite{harrell02}, which was one of the sources of inspiration for this paper. I also thank Dr. Fumiya Okazaki for helpful discussions on Riemannian geometry. The author's work is partially supported by JSPS KAKENHI JP25KJ0120.

\bibliographystyle{abbrv}
\bibliography{ref}

\end{document}